\documentclass{article}

\usepackage{arxiv}

\usepackage[utf8]{inputenc}
\usepackage[T1]{fontenc}
\usepackage{hyperref}
\usepackage{url}
\usepackage{booktabs}
\usepackage{amsfonts}
\usepackage{nicefrac}
\usepackage{microtype}
\usepackage{lipsum}
\usepackage{graphicx}
\usepackage{amsmath,amssymb,amsthm,mathtools}
\usepackage{enumitem}
\usepackage{subcaption}
\usepackage{xcolor}

\theoremstyle{plain}
\newtheorem{theorem}{Theorem}[section]
\newtheorem{lemma}[theorem]{Lemma}
\newtheorem{proposition}[theorem]{Proposition}
\newtheorem{corollary}[theorem]{Corollary}
\newtheorem{problem}[theorem]{Problem}

\theoremstyle{definition}
\newtheorem{definition}[theorem]{Definition}
\newtheorem{remark}[theorem]{Remark}

\usepackage[numbers]{natbib}

\graphicspath{ {./images/} }

\title{A face polynomial invariant for embedded graphs}

\author{
 Hongxia Ma \\
  College of Mathematical Sciences\\
  Xinjiang Normal University\\
  Urumqi 830017, China \\
   \And
 Metrose Metsidik \\
  College of Mathematical Sciences\\
  Xinjiang Normal University\\
  Urumqi 830017, China \\
  \texttt{metrose@xjnu.edu.cn} \\
}

\date{}

\begin{document}
\maketitle

\begin{abstract}
We introduce the face polynomial of a ribbon graph, a generating function for face numbers under the ribbon group action. It addresses the face-count case of the Ellis--Monaghan--Moffatt question on twisted duality. We derive deletion--contraction recursions and prove that this polynomial coincides with the vertex polynomial for \(\langle\delta\rangle\), \(\langle\delta\tau\rangle\), and \(\langle\delta,\tau\rangle\). We characterize monomial cases and obtain formulas and recursions for Euler-genus partial-twuality polynomials. For orientable graphs, we show that the partial-dual face polynomial has both parities whenever an edge exists, so the partial-dual Euler-genus polynomial has only even exponents.
\end{abstract}

\noindent\textbf{Keywords.} twisted duality; ribbon group action; face polynomial; vertex polynomial; Euler-genus partial-twuality polynomial

\noindent\textbf{2000 MR Subject Classification.} primary 05C10, 05C31; secondary 57M15

\section{Introduction}

Chmutov~\cite{chmutov2009generalized} introduced partial duality for embedded graphs by defining a dual operation that acts on individual edges rather than on the whole graph. This notion extends the classical geometric dual and has found considerable applications in both mathematics and physics (see, for example,~\cite{chmutov2009generalized}, \cite{gross2020partial}, \cite{krajewski2011topological} and the references therein). Later, Ellis-Monaghan and Moffatt~\cite{ellis2012twisted} further developed the theory by considering two edge operations: taking the dual with respect to an edge and adding a half-twist to an edge. These operations generate the ribbon group action, which is isomorphic to \(S_3^{|E(G)|}\) and contains geometric duality and Petrie duality as special cases. This group action supplies a powerful framework for analyzing structural features of graph polynomials, especially through the generalized transition polynomial.

Two embedded graphs are called twisted duals if they belong to the same orbit under the ribbon group action. Twisted duality provides a unified perspective for studying several topological graph polynomials. For instance, it allows the Penrose polynomial to be recovered from the transition polynomial~\cite{ellis2012twisted}, and it clarifies the relationship between the Bollobás--Riordan polynomial and the transition polynomial under appropriate restrictions~\cite{ellis2011recipe}. These developments indicate that twisted duality has broad potential for deepening the understanding of various graph-theoretic problems.

Ellis-Monaghan and Moffatt~\cite{ellis2012twisted} posed the following basic question:

\begin{problem}\label{prob:1}
How do invariants of an embedded graph \(G\), such as orientability, genus, degree sequence, or number of vertices, vary within an orbit of the ribbon group action? Is it possible to determine ranges for any of these invariants?
\end{problem}

In connection with the genus case of Problem~\ref{prob:1}, Gross, Mansour, and Tucker~\cite{gross2021partial} introduced the partial-dual polynomial of a ribbon graph, a generating function that enumerates all partial duals according to Euler genus. This polynomial has been widely studied~\cite{chmutov2021conjecture}, \cite{gross2020partial}, \cite{gross2021partial}, \cite{yan2021counterexamples}, \cite{yan2022partial} and has also been extended to delta-matroids~\cite{yan2022twist}, \cite{yan2024partial}, \cite{yuschak2024delta}. For the vertex-count case of Problem~\ref{prob:1}, Yan et al.~\cite{YDM2026Vertex} introduced the vertex polynomial. Their Definition~3.1 indexes the sum by orbit elements, while the paper also gives state-sum characterizations indexed by edge subsets. These conventions are not equivalent in general: distinct edge subsets may parametrize the same ribbon graph in an orbit. We therefore adopt the state-sum convention, counting each parametrizing state with multiplicity, which is the convention compatible with transition-polynomial-type generating functions.

In this paper, we introduce the face polynomial analogously, but recording face numbers under the ribbon group action. While~\cite{YDM2026Vertex} converts vertex-count variation into face-count variation via boundary-component identities, we study face-count variation directly. The face polynomial coincides with the vertex polynomial for \(\langle\delta\rangle\), \(\langle\delta\tau\rangle\), and \(\langle\delta,\tau\rangle\). The case \(\langle\tau\delta\tau\rangle\) reduces immediately, while the case \(\langle\tau\rangle\) requires separate treatment. We also recast the orbit-based arguments of~\cite{YDM2026Vertex} in the state-sum framework and give full proofs, thereby addressing the face-count case of Problem~\ref{prob:1} and establishing the analogous identities and recursions.

The main contributions of the paper are as follows:
\begin{enumerate}[label=(\arabic*)]
\item we introduce the face polynomial for the subgroup actions generated by partial duality and partial Petrial;
\item we derive simplified state-sum formulas and deletion--contraction recursions for the face polynomial;
\item we compare the face polynomial with the vertex polynomial and the transition polynomial, and identify the subgroup actions for which the face and vertex polynomials coincide;
\item we apply the face polynomial to obtain explicit expressions, recursions, and monomial criteria for Euler-genus partial-twuality polynomials;
\item in the orientable case, we analyze the parity structure of the support of the partial-dual face polynomial and prove that both parities occur in its support whenever the graph has at least one edge; consequently, the partial-dual Euler-genus polynomial has only even exponents.
\end{enumerate}

The remainder of this paper is organized as follows. In Section 2, we recall the necessary background on ribbon graphs, the partial dual and partial Petrial operations, and the ribbon group action. Section 3 introduces the face polynomial for various subgroups of the ribbon group, establishes simplified formulas via face-count identities, and derives recursive deletion--contraction relations. In Section 4, we investigate basic properties of the face polynomial, including its evaluation at \(x=1\), invariance under the corresponding orbit, interpolation, and monomial characterizations; we also compare it with the vertex polynomial and the transition polynomial. Section 5 applies the face polynomial to obtain explicit expressions and recursions for the Euler-genus partial-twuality polynomials. Finally, Section 6 focuses on the orientable case, where we analyze the parity structure of the support of the partial-dual face polynomial and derive related results for the partial-dual Euler-genus polynomial.

\section{Background on ribbon graphs and twisted duality}

Ribbon graphs are the standard model for cellularly embedded graphs (see~\cite{ellis2013graphs}) and are the main model used in this paper for graphs embedded in surfaces. We begin by recalling the standard definitions and known results needed throughout.

\begin{definition}[\cite{bollobas2002polynomial}]\label{def:ribbon}
A \emph{ribbon graph} \(G=(V(G),E(G))\) is a compact surface with boundary (possibly non-orientable), formed as the union of two families of topological discs: a set \(V(G)\) of vertex-discs and a set \(E(G)\) of edge-discs. The following conditions are required.
\begin{enumerate}[label=(\arabic*)]
\item Vertex-discs and edge-discs meet only in disjoint line segments;
\item Each such segment lies on the boundary of exactly one vertex-disc and exactly one edge-disc;
\item Every edge-disc contains exactly two such segments.
\end{enumerate}
\end{definition}

We first recall Chmutov's partial dual, which dualizes an embedded graph with respect to a chosen edge subset and is one of the two generators of the ribbon group action.

\begin{definition}[\cite{chmutov2009generalized}]\label{def:partialdual}
Let \(G\) be a ribbon graph and \(A\subseteq E(G)\). The \emph{partial dual} \(G^{\delta(A)}\) is constructed as follows: attach a disc to each boundary component of the spanning ribbon subgraph \((V(G),A)\); these discs become the vertex-discs of \(G^{\delta(A)}\). The interiors of the original vertex-discs of \(G\) are then removed, while all edge-discs remain unchanged.
\end{definition}

The second edge operation, the partial Petrial, twists a chosen edge subset and, together with the partial dual, generates the ribbon group action.

\begin{definition}[\cite{wilson1979operators}]\label{def:petrial}
Let \(G\) be a ribbon graph and \(A\subseteq E(G)\). The \emph{partial Petrial} \(G^{\tau(A)}\) is obtained from \(G\) by inserting a half-twist in every edge belonging to \(A\).
\end{definition}

\begin{figure}[htbp]
  \centering
  \includegraphics[width=0.5\textwidth]{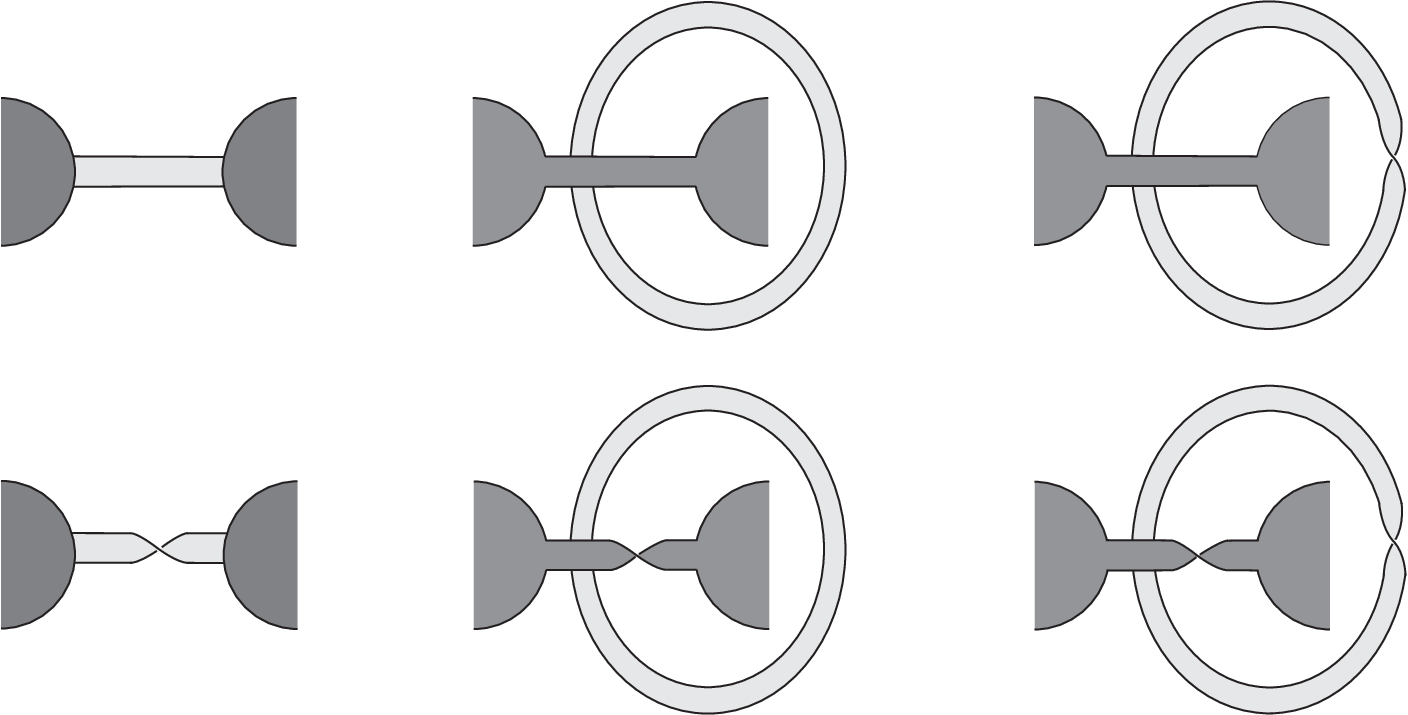}
  \begin{picture}(0,0)
    \put(-175,93){\(\delta\)}
    \put(-178,85){\(\leftrightarrow\)}
    \put(-88,93){\(\tau\)}
    \put(-90,85){\(\leftrightarrow\)}
\put(-175,33){\(\delta\)}
\put(-178,25){\(\leftrightarrow\)}
\put(-86,33){\(\tau\)}
\put(-88,25){\(\leftrightarrow\)}
\put(-210,56){\(\tau\)}
\put(-218,56){\(\updownarrow\)}
\put(-65,55){\(\delta\)}
\put(-57,56){\(\updownarrow\)}
  \end{picture}
  \caption{The effect of \(\delta\) and \(\tau\) on a single edge.}
  \label{fig:edge}
\end{figure}

Figure~\ref{fig:edge} illustrates how \(\delta\) and \(\tau\) act on a single edge. As shown in~\cite{ellis2012twisted}, the partial dual \(\delta\) and the partial Petrial \(\tau\) generate an action of the symmetric group \(S_3\) on ribbon graphs, with presentation
\[
\mathfrak{B} := \langle \delta, \tau \mid \delta^2=\tau^2=(\delta\tau)^3=\mathbf{1} \rangle \cong S_3.
\]

Successive applications to a ribbon graph \(G\) are abbreviated by
\[
G^{\pi(A)\rho(B)}:=\bigl(G^{\pi(A)}\bigr)^{\rho(B)}
\qquad(A,B\subseteq E(G),\ \pi,\rho\in\mathfrak{B}).
\]

Applying \(\delta\) and \(\tau\) globally recovers the geometric and Petrie duals, whereas restricting them to an edge subset yields their partial counterparts. Iterating these operations produces further classes of ribbon graphs, which are summarized in Table~\ref{tab:twisted} and formalized in Definition~\ref{def:orbits}.

\begin{table}[htbp]
\centering
\caption{Ribbon graph classes arising from \(\delta\) and \(\tau\).}
\label{tab:twisted}
\begin{tabular}{|c|c|c|c|}
\hline
Operation(s) & Order & Action on \(E(G)\) & Action on a subset of \(E(G)\) \\
\hline
\(\delta\) & 2 & geometric dual & partial dual \\
\(\tau\) & 2 & Petrie dual & partial Petrial \\
\(\tau \delta \tau\) & 2 & Wilson dual & partial Wilsonial \\
\(\delta \tau\) or \(\tau \delta\) & 3 & triality & partial triality \\
\(\delta\) and \(\tau\) & 6 & a direct derivative & twisted dual \\
\hline
\end{tabular}
\end{table}

The notation \(\biguplus_{i=1}^n A_i=E(G)\) means that \(A_1,\dots,A_n\) are pairwise disjoint and their union is \(E(G)\); some of the \(A_i\) may be empty.

\begin{definition}[\cite{ellis2012twisted}]\label{def:orbits}
Let \(G\) be a ribbon graph. For \(\pi\in\{\delta,\tau,\tau\delta\tau\}\), set
\[
\mathrm{Orb}_{\langle\pi\rangle}(G):=\{ G^{\pi(A)} \mid A\subseteq E(G)\}.
\]
For \(\pi\in\{\delta\tau,\tau\delta\}\), set
\[
\mathrm{Orb}_{\langle\pi\rangle}(G):=\{ G^{\mathbf{1} (A_1)\tau\delta(A_2)\delta\tau(A_3)} \mid \biguplus_{i=1}^3 A_i=E(G)\}.
\]
Finally, define
\[
\mathrm{Orb}_{\langle\delta,\tau\rangle}(G):=\{ G^{\mathbf{1} (A_1)\delta(A_2)\tau(A_3)\tau\delta(A_4)\delta\tau(A_5)\tau\delta\tau(A_6)} \mid \biguplus_{i=1}^6 A_i=E(G)\}.
\]
A ribbon graph \(H\) is called a \emph{twisted dual} of \(G\) if \(H\in\mathrm{Orb}_{\langle\delta,\tau\rangle}(G)\).
\end{definition}

The following basic properties of twisted duality will be used repeatedly; in particular, they allow us to reduce many computations to single edges and to connected components.

\begin{proposition}[Prop.~2.5 in \cite{ellis2012twisted}]\label{pro:twisted-duals}
Let \(G\) be a ribbon graph, \(A,B\subseteq E(G)\), and \(\pi,\rho\in\mathfrak{B}\). Then:
\begin{enumerate}[label=(\arabic*)]
\item If \(A\cap B=\emptyset\), then \(G^{\pi(A)\rho(B)}=G^{\rho(B)\pi(A)}\);
\item If \(e\in A\), then \(G^{\pi(A)}=\bigl(G^{\pi(e)}\bigr)^{\pi(A\setminus\{e\})}\).
\end{enumerate}
\end{proposition}

\section{The face polynomial of ribbon graphs}

To capture the effect of the ribbon group action on face counts, we introduce the following generating function. Here \(f(G)\) denotes the number of boundary components of a ribbon graph \(G\).

\begin{definition}\label{def:facepoly}
For \(\bullet \in \{\langle \delta \rangle, \langle \tau \rangle, \langle \delta\tau \rangle, \langle \tau\delta \rangle, \langle \tau\delta\tau \rangle, \langle \delta,\tau \rangle\}\), the face polynomial for the action \(\bullet\) of a ribbon graph \(G\) is defined by
\begin{align*}
FP_{\langle\pi\rangle}(G,x) 
&= \sum_{A\subseteq E(G)} x^{f(G^{\pi(A)})}, \quad \text{for } \pi \in \{\delta, \tau, \tau\delta\tau\};\\
FP_{\langle\pi\rangle}(G,x) 
&= \sum_{\biguplus_{i=1}^3 A_i = E(G)} x^{f\big(G^{\mathbf{1}(A_1)\tau\delta(A_2)\delta\tau(A_3)}\big)}, 
\quad \text{for } \pi \in \{\delta\tau, \tau\delta\};\\
FP_{\langle\delta,\tau\rangle}(G,x) 
&= \sum_{\biguplus_{i=1}^6 A_i = E(G)} x^{f\big(G^{\mathbf{1}(A_1)\delta(A_2)\tau(A_3)\tau\delta(A_4)\delta\tau(A_5)\tau\delta\tau(A_6)}\big)}.
\end{align*}
\end{definition}

Two immediate reductions follow from the definition and simplify the subsequent discussion.

\begin{remark}\label{rem:reductions}
Let \(G\) be a ribbon graph. The Wilsonial operation \(\tau\delta\tau\) preserves the face count, that is,
\[
f(G^{\tau\delta\tau(A)}) = f(G)
\]
for every subset \(A \subseteq E(G)\). This follows from the local description in Figure~\ref{fig:edge} together with the edge-by-edge compatibility in Proposition~\ref{pro:twisted-duals}(2); see also~\cite{ellis2012twisted}. Consequently,
\[
FP_{\langle\tau\delta\tau\rangle}(G,x) = 2^{e(G)} x^{f(G)}.
\]
Furthermore, \(\langle\delta\tau\rangle=\langle\tau\delta\rangle\) implies
\[
FP_{\langle\delta\tau\rangle}(G,x) = FP_{\langle\tau\delta\rangle}(G,x).
\]
Therefore, it suffices to study the face polynomial for the cases
\[
\bullet \in \{\langle \delta \rangle, \langle \tau\rangle, \langle \delta\tau \rangle, \langle \delta,\tau \rangle\}
\]
in the subsequent analyses.
\end{remark}

To translate face numbers under the group action into those of spanning subgraphs, we need the following identities. Here \(G-A\) denotes the ribbon subgraph obtained by deleting the edges in \(A\) while retaining all vertices.

\begin{lemma}\label{lem:vertexface}
Let \(G\) be a ribbon graph and \(A \subseteq E(G)\). Then
\[
f(G^{\delta(A)}) = f(G^{\delta\tau(A)}) = f(G - A)
\]
and
\[
f(G^{\tau(A)}) = f(G^{\tau\delta(A)})=f(G^{\tau\delta\tau(A)}- A)=f(G^{\delta\tau(A)}- A).
\]
\end{lemma}
\begin{proof}
By Proposition~\ref{pro:twisted-duals}(2), it remains only to verify the case \(A=\{e\}\), which is readily verified from Figure~\ref{fig:face-equiv}.
\begin{figure}[htbp]
  \centering
  \includegraphics[width=0.5\textwidth]{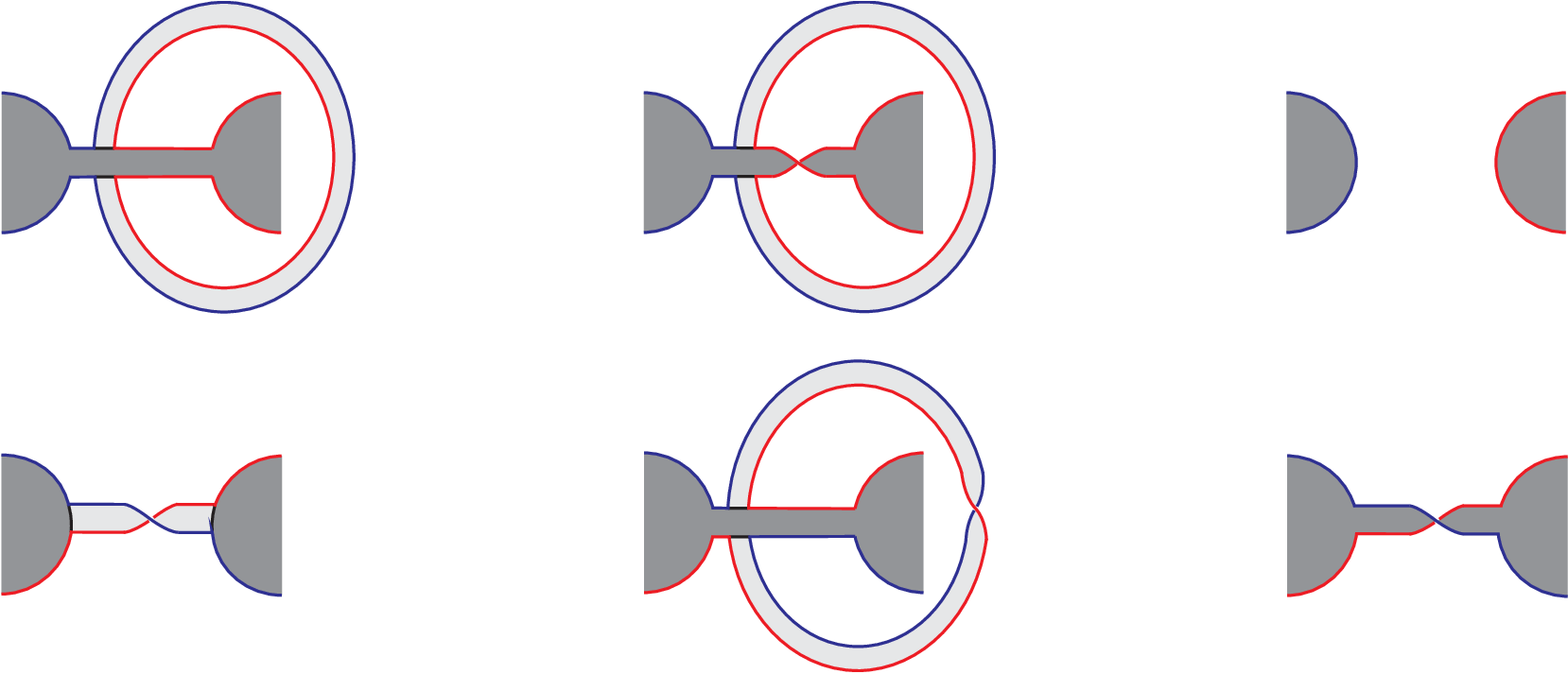}
  \begin{picture}(0,0)
    \put(-240,92){\({G^{\delta(e)}}\)}
    \put(-145,95){\({G^{\delta\tau(e)}}\)}
    \put(-43,92){\({G-e}\)}
    \put(-226,37){\({G^{\tau(e)}}\)}
    \put(-148,40){\({G^{\tau\delta(e)}}\)}
    \put(-43,0){\({G^{\tau\delta\tau(e)}-e}\)}
    \put(-43,37){\({G^{\delta\tau(e)}-e}\)}
    \put(-253,74){\({f\bigg(}\)}
    \put(-180,74){\({\bigg)=f\bigg(}\)}
    \put(-85,74){\({\bigg)=f\bigg(}\)}
    \put(0,74){\({\bigg)}\)}
    \put(-253,19){\({f\bigg(}\)}
    \put(-187,19){\({\bigg)\ =\ f\bigg(}\)}
    \put(-82,19){\({\bigg)=f\bigg(}\)}
    \put(0,19){\({\bigg)}\)}
  \end{picture}
  \caption{Effect of the operations on faces.}
  \label{fig:face-equiv}
\end{figure}
\end{proof}

The next lemma provides a useful normal form for twisted duals and will be used to rewrite the state sums below. Recall that, for \(A\subseteq E(G)\), \(G/A\) denotes the ribbon graph obtained by contracting all edges in \(A\) in the sense of~\cite{ellis2013graphs}, and \(G/e:=G/\{e\}\) for \(e\in E(G)\).

\begin{lemma}\label{lem:contractionedge-face}
Let \(G\) be a ribbon graph and \(A,B \subseteq E(G)\). Then
\[
f((G/A)^{\delta(B)}) = f(G^{\delta(B\setminus A)}).
\]
\end{lemma}
\begin{proof}
\[
\begin{aligned}
f\bigl((G / A)^{\delta(B)}\bigr) 
&= f\bigl((G^{\delta(A)} - A) - B\bigr) = f\bigl(G^{\delta(A)} - (A \cup B)\bigr)\\& = f\bigl((G^{\delta(A)})^{\delta(A \cup B)}\bigr) = f\bigl(((G^{\delta(A)})^{\delta(A)})^{\delta(B\setminus A)}\bigr) = f(G^{\delta(B\setminus A)}).
\end{aligned}
\]
\end{proof}

Combining the reductions in Remark~\ref{rem:reductions} with Lemmas~\ref{lem:vertexface} and~\ref{lem:contractionedge-face}, we obtain the following face-count expressions for the face polynomial.

\begin{theorem}\label{thm:faceequiv}
Let \(G\) be a ribbon graph. Then the following hold:
\begin{enumerate}[label=(\arabic*)]
\item \(FP_{\langle\delta\rangle}(G,x) = \sum_{A\subseteq E(G)} x^{f(G- A)}\);
\item \(FP_{\langle\tau\rangle}(G,x) = \sum_{A\subseteq E(G)} x^{f(G^{\tau\delta\tau(A)}- A)} = \sum_{A\subseteq E(G)} x^{f(G^{\delta\tau(A)}- A)}\);
\item \(FP_{\langle\delta\tau\rangle}(G,x) = \sum_{\biguplus_{i=1}^3 A_i = E(G)} x^{f(G^{\tau(A_2)}- A_3)}\);
\item \(FP_{\langle\delta,\tau\rangle}(G,x) = \sum_{\biguplus_{i=1}^6 A_i = E(G)} x^{f(G^{\tau(A_3\cup A_4)}- (A_2\cup A_5))}\).
\end{enumerate}
\end{theorem}
\begin{proof}
(1)--(2) For \(\pi \in \{\delta, \tau\}\), by Lemma~\ref{lem:vertexface}, the face count \(f(G^{\pi(A)})\) corresponds to \(f(G- A)\) for \(\pi=\delta\) and to \(f(G^{\tau\delta\tau(A)}- A)\) or \(f(G^{\delta\tau(A)}- A)\) for \(\pi=\tau\). Thus,
\begin{align*}
FP_{\langle\delta\rangle}(G,x)&=\sum_{A\subseteq E(G)} x^{f(G- A)},\\
FP_{\langle\tau\rangle}(G,x)&=\sum_{A\subseteq E(G)} x^{f(G^{\tau\delta\tau(A)}- A)}
= \sum_{A\subseteq E(G)} x^{f(G^{\delta\tau(A)}- A)}.
\end{align*}

(3) 
\[
f(G^{1(A_1)\tau\delta(A_2)\delta\tau(A_3)}) = f(G^{\tau\delta(A_2)}-A_3)  = f(G^{\tau(A_2)}-A_3),
\]

(4) 
\[
\begin{aligned}
f(G^{1(A_1)\delta(A_2)\tau(A_3)\tau\delta(A_4)\delta\tau(A_5)\tau\delta\tau(A_6)}) 
&= f(G^{\tau(A_3\cup A_4)\delta(A_2\cup A_5)}) = f(G^{\tau(A_3\cup A_4)}-(A_2\cup A_5)).
\end{aligned}
\]
\end{proof}

As a consequence of Theorem~\ref{thm:faceequiv}, the face polynomial for the full ribbon group action can be reduced to that for the \(\langle\delta\tau\rangle\)-action.

\begin{proposition}\label{prop:connections}
Let \(G\) be a ribbon graph. Then
\begin{align}
FP_{\langle\delta,\tau\rangle}(G,x) &= 2^{e(G)} FP_{\langle\delta\tau\rangle}(G,x).
\end{align}
\end{proposition}

\begin{proof}
From Theorem~\ref{thm:faceequiv}, we relabel the six disjoint sets \(A_1,\dots,A_6\) as
\[
B_1:=A_1\cup A_6,\quad B_2:=A_3\cup A_4,\quad B_3:=A_2\cup A_5.
\]
Under the original partition, every edge lying in \(B_i\) (\(i=1,2,3\)) corresponds to two independent choices. Thus each \(B_i\) contributes a factor of \(2^{|B_i|}\). Therefore,
\begin{align*}
FP_{\langle\delta,\tau\rangle}(G,x) &= \sum_{\biguplus_{i=1}^3 B_i = E(G)} 2^{|B_1|+|B_2|+|B_3|} x^{f(G^{\tau(B_2)}- B_3)}\\
&= 2^{e(G)} \sum_{\biguplus_{i=1}^3 B_i = E(G)} x^{f(G^{\tau(B_2)}- B_3)} = 2^{e(G)} FP_{\langle\delta\tau\rangle}(G,x).
\end{align*}
\end{proof}

Consequently, we only need to investigate the face polynomial for
\[
\bullet \in \{\langle \delta \rangle, \langle \tau\rangle, \langle \delta\tau \rangle\}.
\]

We now show that the face polynomial obeys deletion--contraction type recursions, reflecting the edge operations that generate the ribbon group action.

\begin{theorem}\label{thm:recursions}
Let \(G\) be a ribbon graph and \(e \in E(G)\). Then
\begin{enumerate}[label=(\arabic*)]
\item \(FP_{\langle\delta\rangle}(G,x) = FP_{\langle\delta\rangle}(G-e,x) +FP_{\langle\delta\rangle}(G/e,x)\);
\item \(FP_{\langle\tau\rangle}(G,x) = FP_{\langle\tau\rangle}(G/e,x) + FP_{\langle\tau\rangle}(G^{\tau(e)}/e,x)\);
\item \(FP_{\langle\delta\tau\rangle}(G,x) = FP_{\langle\delta\tau\rangle}(G- e,x) + FP_{\langle\delta\tau\rangle}(G/e,x) + FP_{\langle\delta\tau\rangle}(G^{\tau(e)}/e,x)\).
\end{enumerate}
\end{theorem}
\begin{proof}
Partitioning edge subsets by \(e\):
\begin{align*}
FP_{\langle\pi\rangle}(G,x) &= \sum_{A\subseteq E(G)} x^{f(G^{\pi(A)})} 
= \sum_{\substack{A\subseteq E(G)\\ e\in A}} x^{f(G^{\pi(A)})} + \sum_{\substack{A\subseteq E(G)\\ e\notin A}} x^{f(G^{\pi(A)})}\\
&= \sum_{A'\subseteq E(G)\setminus \{e\}} x^{f((G^{\pi(e)})^{\pi(A')})}+\sum_{A\subseteq E(G)\setminus \{e\}} x^{f(G^{\pi(A)})},
\end{align*}
where \(\pi\in\{\delta,\tau\}\).

(1) By Lemma~\ref{lem:vertexface}, for \(A',A\subseteq E(G)\setminus \{e\}\):
\[
\begin{aligned}
f\big((G-e)^{\delta(A')}\big) 
&= f\big((G-e)-A'\big) = f\big(G-(\{e\}\cup A')\big) = f\big(G^{\delta(\{e\}\cup A')}\big) = f\big((G^{\delta(e)})^{\delta(A')}\big).
\end{aligned}
\]
and, by Lemma~\ref{lem:contractionedge-face}, we have
\[
\begin{aligned}
f\big((G/e)^{\delta(A)}\big) = f\big(G^{\delta(A\setminus\{e\})}\big) = f\big(G^{\delta(A)}\big).
\end{aligned}
\]
Hence,
\begin{align*}
FP_{\langle\delta\rangle}(G,x) &= \sum_{A'\subseteq E(G-e)} x^{f\big((G-e)^{\delta(A')}\big)} + \sum_{A\subseteq E(G/e)} x^{f\big((G/e)^{\delta(A)}\big)}= FP_{\langle\delta\rangle}(G-e,x) + FP_{\langle\delta\rangle}(G/e,x).
\end{align*}

(2) 
\[
\begin{aligned}
f\big((G^{\tau(e)}/e)^{\tau(A')}\big) 
&= f\big((G^{\delta\tau(e)}- e)^{\delta\tau(A')}-A'\big) \\
&= f\big(G^{\delta\tau(\{e\}\cup A')}- (\{e\}\cup A')\big) \\
&= f\big(G^{\delta\tau(\{e\}\cup A')\delta\tau(\{e\}\cup A')}\big)\\
&= f\big(G^{\tau\delta(\{e\}\cup A')}\big)\\
&= f\big(G^{\tau(\{e\}\cup A')}\big)\\
&= f\big((G^{\tau(e)})^{\tau(A')}\big),
\end{aligned}
\]
and, by Lemma~\ref{lem:contractionedge-face}, we have
\[
\begin{aligned}
f\big((G/e)^{\tau(A)}\big) = f\big(G^{\tau(A)}/e\big) = f\big(G^{\tau(A)}\big).
\end{aligned}
\]
Hence,
\begin{align*}
FP_{\langle\tau\rangle}(G,x) &= \sum_{A'\subseteq E(G^{\tau(e)}/e)} x^{f\big((G^{\tau(e)}/e)^{\tau(A')}\big)} + \sum_{A\subseteq E(G/e)} x^{f\big((G/e)^{\tau(A)}\big)}\\
&= FP_{\langle\tau\rangle}(G^{\tau(e)}/e,x) + FP_{\langle\tau\rangle}(G/e,x).
\end{align*}

(3) By Theorem~\ref{thm:faceequiv}, we decompose the sum based on the position of edge \(e\):
\begin{align*}
FP_{\langle\delta\tau\rangle}(G,x) &= \sum_{\biguplus_{i=1}^3 A_i = E(G)} x^{f\big(G^{\tau(A_2)}- A_3\big)} \\
&= \sum_{\substack{\biguplus A_i=E(G)\\ e\in A_3}}x^{f\big(G^{\tau(A_2)}- A_3\big)} +\sum_{\substack{\biguplus A_i=E(G)\\ e\in A_1}}x^{f\big(G^{\tau(A_2)}- A_3\big)} + \sum_{\substack{\biguplus A_i=E(G)\\ e\in A_2}}x^{f\big(G^{\tau(A_2)}- A_3\big)}.
\end{align*}

Case 1: \(e\in A_3\).
\[
\begin{aligned}
\sum_{\substack{\biguplus A_i=E(G)\\ e\in A_3}} x^{f\big(G^{\tau(A_2)}- A_3\big)}
&= \sum_{\biguplus A_i = E(G-e)} x^{f\big(G^{\tau(A_2)}- A_3-e\big)} = \sum_{\biguplus A_i = E(G-e)} x^{f\big((G-e)^{\tau(A_2)}\setminus A_3\big)} = FP_{\langle\delta\tau\rangle}(G-e,x).
\end{aligned}
\]

Case 2: \(e\in A_1\).
\[
\begin{aligned}
\sum_{\substack{\biguplus A_i=E(G)\\ e\in A_1}} x^{f\big(G^{\tau(A_2)}- A_3\big)}
&= \sum_{\biguplus A_i = E(G/e)} x^{f\big((G^{\tau(A_2)}-A_3)/e\big)} = \sum_{\biguplus A_i = E(G/e)} x^{f\big((G/e)^{\tau(A_2)}- A_3\big)} = FP_{\langle\delta\tau\rangle}(G/e,x).
\end{aligned}
\]

Case 3: \(e\in A_2\).
\[
\begin{aligned}
\sum_{\substack{\biguplus A_i=E(G)\\ e\in A_2}} x^{f\big(G^{\tau(A_2)}- A_3\big)}
&= \sum_{\biguplus A_i = E(G^{\tau(e)}/e)} x^{f\big((G^{\tau(e)})^{\tau(A_2)}- A_3\big)} \\&= \sum_{\biguplus A_i = E(G^{\tau(e)}/e)} x^{f\big((G^{\tau(e)}/e)^{\tau(A_2)}- A_3\big)} = FP_{\langle\delta\tau\rangle}(G^{\tau(e)}/e,x).
\end{aligned}
\]
\end{proof}

\section{Basic properties of the face polynomial}

For a ribbon graph \(G\), let \(v(G)\) and \(e(G)\) denote its numbers of vertices and edges, respectively. For a polynomial
\(
f(x)=\sum_{i=0}^m a_i x^i,
\)
its support is defined by
\(
\operatorname{Supp}(f)=\{i : a_i \neq 0\},
\)
the set of indices of the non-zero coefficients. We call a polynomial with nonnegative coefficients \emph{interpolating} if its support contains every integer between its minimum and maximum degrees.

We begin with basic evaluations and structural properties that follow directly from the state-sum definition and orbit invariance.

\begin{proposition}\label{prop:evaluations}
Let \(G\) be a connected ribbon graph. Then
\begin{enumerate}[label=(\arabic*)]
\item \(FP_{\langle\delta\rangle}(G,1) = FP_{\langle\tau\rangle}(G,1) = 2^{e(G)}\) and  \(FP_{\langle\delta\tau\rangle}(G,1) = 3^{e(G)}\);
\item For any \(\bullet \in \{\langle \delta\rangle,\langle \tau\rangle,\langle \delta\tau\rangle\}\), \(FP_{\bullet}(G,x) = FP_{\bullet}(H,x)\) holds for all \(H \in \mathrm{Orb}_{\bullet}(G)\);
\item For \(\bullet \in \{\langle \delta\rangle,\langle \tau\rangle,\langle \delta\tau\rangle\}\), the face polynomial \(FP_{\bullet}(G,x)\) is interpolating with minimum degree \(1\).
\end{enumerate}
\end{proposition}
\begin{proof}
(1) Evaluating \(FP_{\langle\delta\rangle}(G,1)\) and \(FP_{\langle\tau\rangle}(G,1)\) yields the total number of edge subsets of the graph \(G\), namely \(2^{e(G)}\). On the other hand, \(FP_{\langle\delta\tau\rangle}(G,1)\) enumerates all partial trialities of \(G\) and equals \(3^{e(G)}\), since each edge admits exactly three choices: the identity, \(\delta\tau\), and \(\tau\delta\).

(2) Let \(H \in \mathrm{Orb}_{\langle \delta\rangle}(G)\). Then \(H = G^{\delta(A)}\) for some \(A \subseteq E(G)\). Thus
\[
G^{\delta(B)} = H^{\delta(A\triangle B)}
\]
for all \(B \subseteq E(G)\). Clearly, the map \(B \mapsto A\triangle B\) is a bijection, so \(FP_{\langle \delta\rangle}(G,x) = FP_{\langle \delta\rangle}(H,x)\). Analogously, \(FP_{\langle\tau\rangle}(G,x) = FP_{\langle\tau\rangle}(H,x)\).

Now let \(H \in \mathrm{Orb}_{\langle \delta\tau\rangle}(G)\). Then \(H = G^{\tau\delta(A_1)\delta\tau(A_2)}\) for two disjoint edge subsets \(A_1,A_2 \subseteq E(G)\). For any pair of disjoint edge subsets \(B_1,B_2 \subseteq E(G)\), we have
\[
\begin{aligned}
G^{\tau\delta(B_1)\delta\tau(B_2)}
&= \bigl(H^{\delta\tau(A_1)\tau\delta(A_2)}\bigr)^{\tau\delta(B_1)\delta\tau(B_2)} \\
&= H^{\delta\tau(A_1)\bigl(\tau\delta(A_2)\tau\delta(B_1)\bigr)\delta\tau(B_2)} \\
&= H^{\delta\tau(A_1)\bigl(\tau\delta(A_2\triangle B_1)\delta\tau(A_2\cap B_1)\bigr)\delta\tau(B_2)} \\
&= H^{\delta\tau(A_1)\delta\tau\bigl(B_2\setminus(A_2\triangle B_1)\bigr)\delta\tau(A_2\cap B_1)\tau\delta\bigl((A_2\triangle B_1)\setminus B_2\bigr)} \\
&= H^{\delta\tau(A_1)\delta\tau\bigl(B_2\setminus A_2\bigr)\delta\tau(A_2\cap B_1)\tau\delta\bigl((A_2\setminus B_2)\triangle B_1\bigr)} \\
&= H^{\delta\tau\bigl(A_1\triangle (B_2\setminus A_2)\bigr)\tau\delta\bigl(A_1\cap (B_2\setminus A_2)\bigr)\delta\tau(A_2\cap B_1)\tau\delta\bigl((A_2\setminus B_2)\triangle B_1\bigr)} \\
&= H^{\delta\tau\Bigl(\bigl(A_1\triangle (B_2\setminus A_2)\bigr)\cup (A_2\cap B_1)\Bigr)
\tau\delta\Bigl(\bigl(A_1\cap B_2\bigr)\cup \bigl((A_2\setminus B_2)\triangle B_1\bigr)\Bigr)} \\
&= H^{\tau\delta\Bigl(\bigl((A_2\setminus B_2)\triangle B_1\bigr)\cup (A_1\cap B_2)\Bigr)
\delta\tau\Bigl(\bigl(A_1\triangle (B_2\setminus A_2)\bigr)\cup (A_2\cap B_1)\Bigr)}
\end{aligned}
\]
The mappings \(B_1 \mapsto \bigl((A_2\setminus B_2)\triangle B_1\bigr)\cup (A_1\cap B_2)\) and \(B_2 \mapsto \bigl(A_1\triangle (B_2\setminus A_2)\bigr)\cup (A_2\cap B_1)\) are both bijections. Consequently, \(FP_{\langle \delta\tau\rangle}(G,x) = FP_{\langle \delta\tau\rangle}(H,x)\).

Therefore, for each \(\bullet \in \{\langle\delta\rangle,\langle\tau\rangle,\langle\delta\tau\rangle\}\), the face polynomial \(FP_{\bullet}(G,x)\) depends only on the orbit of \(G\). It therefore takes the same value for every graph \(H\) in \(\mathrm{Orb}_{\bullet}(G)\).

(3) Let \(H \in \mathrm{Orb}_{\bullet}(G)\) satisfy that the maximum degree of \(FP_{\bullet}(G,x)\) is equal to \(f(H)\). The case \(f(H)=1\) is trivial, so we assume \(f(H)\geq 2\). There then exists an edge \(e\) that is incident to two distinct faces. Applying one of the operators \(\delta(e)\), \(\tau(e)\), \(\delta\tau(e)\) merges these two faces into one, which decreases the total number of faces by \(1\). As a result, \(FP_{\bullet}(H,x)\) is an interpolating polynomial of minimum degree \(1\).
\end{proof}

\begin{theorem}[Theorem~3.5 of \cite{YDM2026Vertex}]\label{YDM2026Vertex}
Let \(G\) be a ribbon graph. Then the vertex polynomial of \(G\) is given by:
\begin{enumerate}[label=(\arabic*)]
    \item \(VP_{\langle\delta\rangle}(G, x)=\sum_{A \subseteq E(G)} x^{f\left(G - A^{c}\right)}\);
    \item \(VP_{\langle\delta \tau\rangle}(G, x)=\sum_{\biguplus_{i=1}^{3} A_{i}=E(G)} x^{f\left(G^{\tau\left(A_{3}\right)} - A_{1}\right)}\);
    \item \(VP_{\langle\delta,\tau\rangle}(G,x) = \sum_{\biguplus_{i=1}^6 A_i = E(G)} x^{f(G^{\tau(A_5\cup A_6)}- (A_1\cup A_3))}\);
    \item \(VP_{\langle\tau\delta\tau\rangle}(G,x) = \sum_{A \subseteq E(G)} x^{f\left(G^\times - A^{c}\right)}=VP_{\langle\delta\rangle}(G^\times, x)\),
\end{enumerate}
\noindent where \(A^c=E(G)\setminus A\) and \(G^\times=G^{\tau(E(G))}\).
\end{theorem}

The preceding formulas allow us to compare the face polynomial with the vertex polynomial. By Theorems~\ref{thm:faceequiv} and~\ref{YDM2026Vertex}, we have

\begin{theorem}\label{Face=Vertex}
Let \(G\) be a ribbon graph. Then:
\begin{enumerate}[label=(\arabic*)]
    \item \(FP_{\langle\delta\rangle}(G, x)=VP_{\langle\delta\rangle}(G, x)\);
    \item \(FP_{\langle\delta \tau\rangle}(G, x)=VP_{\langle\delta \tau\rangle}(G, x)\);
    \item \(FP_{\langle\delta,\tau\rangle}(G, x)=VP_{\langle\delta,\tau\rangle}(G, x)\).
\end{enumerate}
\end{theorem}

This theorem shows that, for the subgroup actions considered here, the face and vertex polynomials share the same state-sum structure. Consequently, several properties of the vertex polynomial derived in~\cite{YDM2026Vertex} carry over directly to the face polynomial.

\begin{proposition}\label{prop:monomial}
Let \(G\) be a connected ribbon graph with \(E(G)\neq \emptyset\). Then
\begin{enumerate}[label=(\arabic*)]
\item \(FP_{\langle\delta\rangle}(G,x)\) is a monomial if and only if \(G\) is a bouquet in which every loop is a trivial non-orientable loop;
\item \(FP_{\langle\tau\rangle}(G,x)\) is a monomial if and only if \(G\) is a tree;
\item \(FP_{\langle\delta\tau\rangle}(G,x)\) is never a monomial.
\end{enumerate}
\end{proposition}
\begin{proof}
Items (1) and (3) are immediate consequences of Theorem~\ref{Face=Vertex} and Proposition 4.4 from~\cite{YDM2026Vertex}. We now proceed to prove assertion (2).

For sufficiency, let \(G\) be a tree with \(m\) edges. Then \(G = G^{\tau(A)}\) for any \(A \subseteq E(G)\). A direct computation gives
\[
FP_{\langle\tau\rangle}(G,x) = \sum_{A\subseteq E(G)} x^{f(G^{\tau(A)})}= 2^m x,
\]
which explicitly shows that the polynomial is a monomial.

For necessity, assume that \(FP_{\langle\tau\rangle}(G,x)\) is a monomial. Proposition~\ref{prop:evaluations}(1) and (3) then give \(FP_{\langle\tau\rangle}(G,x)=2^m x\). Suppose, for contradiction, that \(G\) is not a tree. Then \(G\) contains a cycle. By contracting all edges of a spanning tree of \(G\), we obtain a bouquet \(B\). Clearly, \(B\ncong K_1\). If \(B\) has a trivial loop (i.e., a loop not interlaced with any other loop), then \(FP_{\langle\tau\rangle}(B,x)\) is not a monomial. If \(B\) has exactly two interlaced loops that are not interlaced with any other loop, then \(FP_{\langle\tau\rangle}(B,x)\) is again not a monomial. If \(B\) has more than two interlaced loops, we can choose two of these loops that lie at the ends of the interlaced family (the ``end loops'') and contract all the other loops. This yields a bouquet \(B'\) consisting only of the two chosen end loops; in each case---whether the two loops form a trivial 2-cycle, two interlaced loops, or two trivial loops---\(FP_{\langle\tau\rangle}(B',x)\) is not a monomial. By Theorem~\ref{thm:recursions}(2), the deletion--contraction recursion reduces the computation to \(B\) or \(B'\). Since all coefficients are nonnegative, no cancellation can occur. This implies that \(FP_{\langle\tau\rangle}(G,x)\) is not a monomial, a contradiction. Therefore, \(G\) must be a tree.
\end{proof}

\begin{remark}\label{rem:transition}
It is shown in~\cite{YDM2026Vertex} that the vertex polynomial agrees with the topological transition polynomial \(Q(G;(\boldsymbol{\alpha},\boldsymbol{\beta},\boldsymbol{\gamma}), x)\) introduced in~\cite{ellis2012twisted} under certain specializations. By Theorem~\ref{thm:recursions}, the face polynomial admits a similar identification with the same transition polynomial when the indeterminates \(\boldsymbol{\alpha},\boldsymbol{\beta},\boldsymbol{\gamma}\) are restricted to values in \(\{0,1\}\). Moreover, Theorem~\ref{Face=Vertex} shows that the face and vertex polynomials coincide for the subgroup actions \(\langle\delta\rangle\), \(\langle\delta\tau\rangle\), and \(\langle\delta,\tau\rangle\). These observations place the face polynomial within the broader framework of transition polynomials and clarify its relation to the vertex polynomial. The main contribution of this paper is a direct face-count formulation and a set of tools that streamline the study of Euler-genus partial-twuality polynomials; the parity results in Section~6 illustrate this utility. These connections also highlight the interplay among the face polynomial, the vertex polynomial, and the transition polynomial, and provide complementary perspectives on our results.
\end{remark}

We now apply the face polynomial to the Euler-genus partial-twuality polynomials.

\section{Applications to Euler-genus partial-twuality polynomials}

The Euler genus of a ribbon graph \(G\) is defined as  
\begin{equation}\label{eq:eulergenus}
\epsilon(G)=2c(G)+e(G)-f(G)-v(G),
\end{equation}
where \(c(G)\) is the number of connected components. It equals twice the genus if the embedding is orientable, and equals the crosscap number if the embedding is non-orientable. Recall that Gross, Mansour, and Tucker~\cite{gross2021partial} introduced the partial-dual polynomial of a ribbon graph as a generating function that enumerates all partial duals of the ribbon graph according to Euler genus. Since the effect of individual group elements has been examined, it is natural to study the action of the whole subgroup on the Euler genus.

\begin{definition}\label{def:Euler-genus partial-twuality}
For \(\bullet \in \{\langle \delta \rangle, \langle \tau \rangle, \langle \delta\tau \rangle, \langle \tau\delta\tau \rangle, \langle \delta,\tau \rangle\}\), the \emph{Euler-genus partial \(\bullet\) polynomial} of a ribbon graph \(G\) is defined by
\begin{align*}
\partial \epsilon_G^\bullet(x) &= \sum_{A \subseteq E(G)} x^{\epsilon(G^{\bullet(A)})}, \quad \text{for } \bullet\in \{\langle \delta \rangle, \langle \tau \rangle,  \langle \tau\delta\tau \rangle\};\\
\partial \epsilon_G^{\langle\delta\tau\rangle}(x)
&= \sum_{\biguplus_{i=1}^3 A_i = E(G)} x^{\epsilon\big(G^{\mathbf{1}(A_1)\tau\delta(A_2)\delta\tau(A_3)}\big)};\\
\partial \epsilon_G^{\langle\delta,\tau\rangle}(x)
&= \sum_{\biguplus_{i=1}^6 A_i = E(G)} x^{\epsilon\big(G^{\mathbf{1}(A_1)\delta(A_2)\tau(A_3)\tau\delta(A_4)\delta\tau(A_5)\tau\delta\tau(A_6)}\big)}.
\end{align*}
\end{definition}

Using Euler's formula~\eqref{eq:eulergenus}, we can express the Euler-genus partial-twuality polynomials in terms of face and vertex counts.

\begin{proposition}\label{prop:partial-twuality-face-vertex}
For \(\bullet \in \{\langle \delta \rangle, \langle \tau \rangle, \langle \tau\delta\tau \rangle\}\),
\[
\partial \epsilon_G^\bullet(x)
= x^{2c(G)+e(G)}\sum_{A \subseteq E(G)} x^{-f(G^{\bullet(A)})} x^{-v(G^{\bullet(A)})}.
\]
Moreover,
\begin{align*}
\partial \epsilon_G^{\langle\delta\tau\rangle}(x)
&= x^{2c(G)+e(G)}\sum_{\biguplus_{i=1}^3 A_i = E(G)} x^{-f(G^{\tau\delta(A_2)\delta\tau(A_3)})} x^{-v(G^{\tau\delta(A_2)\delta\tau(A_3)})},\\
\partial \epsilon_G^{\langle\delta,\tau\rangle}(x)
&= x^{2c(G)+e(G)}\sum_{\biguplus_{i=1}^6 A_i = E(G)} x^{-f(G^{\delta(A_2)\tau(A_3)\tau\delta(A_4)\delta\tau(A_5)\tau\delta\tau(A_6)})} x^{-v(G^{\delta(A_2)\tau(A_3)\tau\delta(A_4)\delta\tau(A_5)\tau\delta\tau(A_6)})}.
\end{align*}
\end{proposition}

Substituting the face and vertex-count identities from Theorems~\ref{thm:faceequiv} and~\ref{YDM2026Vertex} gives the following explicit formulas.

\begin{proposition}\label{pro:partial-twuality-face-vertex}
Let \(G\) be a ribbon graph. Then the following hold:
\begin{align*}
\partial \epsilon_G^{\langle\delta\rangle}(x)
&= x^{2c(G)+e(G)}\sum_{A\subseteq E(G)} x^{-f(G-A)}x^{-f(G-A^c)};\\
\partial \epsilon_G^{\langle\tau\rangle}(x) 
&= x^{2c(G)+e(G)-v(G)}\sum_{A\subseteq E(G)} x^{-f(G^{\tau(A)})}
= x^{2c(G)+e(G)-v(G)}\sum_{A\subseteq E(G)} x^{-f(G^{\delta\tau(A)}-A)}\\
&= x^{2c(G)+e(G)-v(G)}FP_{\langle\tau\rangle}(G,x^{-1});\\
\partial \epsilon_G^{\langle\delta\tau\rangle}(x) 
&= x^{2c(G)+e(G)}\sum_{\biguplus_{i=1}^3 A_i = E(G)} x^{-f(G^{\tau(A_2)}-A_3)}x^{-f(G^{\tau(A_3)}-A_1)};\\
\partial \epsilon_G^{\langle\tau\delta\tau\rangle}(x) 
&= x^{2c(G)+e(G)-f(G)}\sum_{A\subseteq E(G)} x^{-f(G^\times-A^c)}= x^{2c(G)+e(G)-f(G)}VP_{\langle\delta\rangle}(G^\times,x^{-1});\\
\partial \epsilon_G^{\langle\delta,\tau\rangle}(x) 
&= x^{2c(G)+e(G)}\sum_{\biguplus_{i=1}^6 A_i = E(G)} x^{-f(G^{\tau(A_3\cup A_4)}-(A_2\cup A_5))}x^{-f(G^{\tau(A_5\cup A_6)}-(A_1\cup A_3))}.
\end{align*}
\end{proposition}

The recursions for the face polynomial transfer directly to the Euler-genus partial-twuality polynomials. Using Theorems~\ref{thm:recursions}, \ref{YDM2026Vertex}, \ref{Face=Vertex}, and Proposition~\ref{pro:partial-twuality-face-vertex}, we obtain the following recursions.

\begin{proposition}\label{pro:recursions}
Let \(G\) be a ribbon graph and \(e \in E(G)\). Then
\begin{enumerate}[label=(\arabic*)]
\item \(\partial \epsilon_G^{\langle\tau\rangle}(x)= \partial \epsilon_{G/e}^{\langle\tau\rangle}(x) + \partial \epsilon_{G^{\tau(e)}/e}^{\langle\tau\rangle}(x)\);
\item \(\partial \epsilon_G^{\langle\tau\delta\tau\rangle}(x)=\partial \epsilon_{G-e}^{\langle\tau\delta\tau\rangle}(x) + \partial \epsilon_{G^{\tau(e)}/e}^{\langle\tau\delta\tau\rangle}(x)\).
\end{enumerate}
\end{proposition}

Using Theorems~\ref{YDM2026Vertex} and~\ref{Face=Vertex} together with Proposition~\ref{prop:monomial}, we translate the monomial characterizations of the face polynomial into the following monomial criteria for Euler-genus partial-twuality polynomials.

\begin{proposition}\label{pro:partial-twuality-monomial}
Let \(G\) be a ribbon graph. Then the following hold:
\begin{enumerate}[label=(\arabic*)]
\item \(\partial \epsilon_G^{\langle\delta\rangle}(x)\) is a monomial if and only if \(G\) is a bouquet in which every loop is a trivial non-orientable loop;
\item \(\partial \epsilon_G^{\langle\tau\rangle}(x)\) is a monomial if and only if \(G\) is a tree;
\item \(\partial \epsilon_G^{\langle\tau\delta\tau\rangle}(x)\) is a monomial if and only if \(G\) is a bouquet in which every loop is a trivial orientable loop;
\item \(\partial \epsilon_G^{\bullet}(x)\) is never a monomial for \(\bullet \in \{\langle \delta\tau \rangle, \langle \tau\delta \rangle, \langle \delta,\tau \rangle\}\).
\end{enumerate}
\end{proposition}

\section{Supports in the orientable case}

A ribbon graph is \emph{orientable} if its underlying surface is orientable, and non-orientable otherwise.

The face polynomial \(FP_{\langle \delta \rangle}(G,x)\) is known to be interpolating. However, the arithmetic structure of its support remains unexplored. For orientable ribbon graphs, Euler's formula imposes a parity constraint on the face count of any subgraph whose number of edges has a prescribed parity.

We first record a parity constraint for orientable ribbon graphs; it is the key tool for the support analysis in this section.

\begin{lemma}\label{lem:face-parity}
Let \(G\) be an orientable ribbon graph. Then for any \(A_1,A_2\subseteq E(G)\),
\[
f(G-A_1)\equiv f(G-A_2)\pmod 2
\quad\Longleftrightarrow\quad
|A_1|\equiv |A_2|\pmod 2.
\]
\end{lemma}

\begin{proof}
Recall Euler's formula for the Euler genus \(\epsilon(G)\):
\[
f(G)=e(G)-v(G)+2c(G)-\epsilon(G).
\]
If \(G\) is orientable, then \(\epsilon(G)=2g(G)\) is even. Since \(2c(G)\) is also even, reducing this identity modulo \(2\) yields
\[
f(G)\equiv e(G)-v(G)\pmod 2. \tag{2}
\]

Now let \(A\subseteq E(G)\). Deleting the edges in \(A\) preserves the vertex set and orientability, so
\[
v(G-A)=v(G),\qquad e(G-A)=e(G)-|A|,\qquad \text{and } G-A \text{ is orientable}.
\]
Applying (2) to \(G-A\), we obtain
\[
f(G-A)\equiv e(G-A)-v(G-A)
= e(G)-v(G)-|A| \pmod 2.
\]
Since \(e(G)-v(G)\) is a fixed constant independent of \(A\), the parity of \(f(G-A)\) depends only on the parity of \(|A|\). Therefore, for any \(A_1,A_2\subseteq E(G)\),
\[
f(G-A_1)\equiv f(G-A_2)\pmod 2
\]
holds if and only if \(|A_1|\equiv |A_2|\pmod 2\), as claimed.
\end{proof}

This parity lemma has the following immediate consequence.

\begin{corollary}
If \(e(G)>0\), then \(FP_{\langle\delta\rangle}(G,x)\) has both even and odd exponents in its support.
\end{corollary}

\begin{proof}
Choose any \(A\subseteq E(G)\) and an edge \(e\in E(G)\). Since \(|A\triangle\{e\}| = |A|\pm 1\), the congruence in Lemma~\ref{lem:face-parity} gives
\[
f(G-(A\triangle\{e\})) \equiv f(G-A)+1 \pmod 2.
\]
Thus both parities occur. For the case \(E(G)=\varnothing\), only one subset \(A=\varnothing\) exists, giving a monomial \(x^{v(G)}\).
\end{proof}

Since, for \(n \ge 1\), there is a one-to-one correspondence between the even-sized subsets and the odd-sized subsets of an \(n\)-element set, Lemma~\ref{lem:face-parity} implies the following proposition.

\begin{proposition}
Let \(G\) be an orientable ribbon graph with \(e(G)>0\). Then
\[
\sum_{\substack{i \in \operatorname{Supp}\big(FP_{\langle \delta \rangle}(G,x)\big) \\ i \text{ even}}} a_i
=
\sum_{\substack{i \in \operatorname{Supp}\big(FP_{\langle \delta \rangle}(G,x)\big) \\ i \text{ odd}}} a_i
=
2^{e(G)-1},
\]
where \(a_i\) denotes the coefficient of \(x^i\) in \(FP_{\langle \delta \rangle}(G,x)\).
\end{proposition}

\begin{proposition}
Let \(G\) be an orientable ribbon graph. Then
\[
\operatorname{Supp}\bigl(\partial \epsilon_G^{\langle\delta\rangle}(x)\bigr) \subseteq 2\mathbb{N}.
\]
\end{proposition}

\begin{proof}
Let \(i\in \operatorname{Supp}\bigl(\partial \epsilon_G^{\langle\delta\rangle}(x)\bigr)\). By Proposition~\ref{pro:partial-twuality-face-vertex}, there exists an edge subset \(A\subseteq E(G)\) such that
\[
i=2c(G)+e(G)-f(G-A)-f(G-A^c).
\]
Notice that
\[
|A|\equiv |A^c| \pmod 2
\]
if \(e(G)\) is even, and
\[
|A|\equiv |A^c|+1 \pmod 2
\]
if \(e(G)\) is odd. Then, by Lemma~\ref{lem:face-parity}, we have
\[
f(G-A)\equiv f(G-A^c) \pmod 2
\]
if \(e(G)\) is even, and
\[
f(G-A)\equiv f(G-A^c)+1 \pmod 2
\]
if \(e(G)\) is odd. Hence \(e(G)-f(G-A)-f(G-A^c)\) is always even. Therefore \(i\) is even, and the proof is complete.
\end{proof}

Recall that the face polynomial is defined for several subgroup actions: \(\langle\delta\rangle\), \(\langle\tau\rangle\), and \(\langle\delta\tau\rangle\). Their supports are naturally related.

\begin{proposition}
For every ribbon graph \(G\),
\[
\operatorname{Supp}\big(FP_{\langle \delta \rangle}(G,x)\big)\cup\operatorname{Supp}\big(FP_{\langle \tau \rangle}(G,x)\big)
\;\subseteq\;
\operatorname{Supp}\big(FP_{\langle \delta \tau \rangle}(G,x)\big).
\]
\end{proposition}

\begin{proof}
By Theorem~\ref{thm:faceequiv}, every term of \(FP_{\langle \delta \rangle}\) corresponds to a subset \(A\subseteq E(G)\) with face count \(f(G-A)\). In the tripartition \((A_1,A_2,A_3)\) defining \(FP_{\langle \delta \tau \rangle}\), setting \(A_2=\varnothing\) and \(A_3=A\) gives exactly the same face count. Hence each \(\delta\)-exponent appears in the \(\delta\tau\)-polynomial. For \(\langle \tau \rangle\), Lemma~\ref{lem:vertexface} gives
\[
f(G^{\tau(A)})=f(G^{\tau\delta\tau(A)}-A),
\]
and the same tripartition with \(A_2=A\) and \(A_3=\varnothing\) yields the corresponding exponent in \(FP_{\langle\delta\tau\rangle}\).
\end{proof}

\section*{Declaration of competing interest}
The authors declare no competing interests.

\section*{Acknowledgements}
This work is supported by National Natural Science Foundation of China (Grant Number: 12661073) and Natural Science Foundation of Xinjiang (Grant Number: 2024D01A89).

\section*{Data availability}
No data was used for the research described in the article.


\begin{thebibliography}{99}

\bibitem{bollobas2002polynomial} Bollobás, B., Riordan, O. A Polynomial of Graphs on Surfaces. Math. Ann., 2002, 323: 81--96
\bibitem{chmutov2009generalized} Chmutov, S. Generalized Duality for Graphs on Surfaces and the Signed Bollobás--Riordan Polynomial. J. Comb. Theory, Ser. B, 2009, 99: 617--638
\bibitem{chmutov2021conjecture} Chmutov, S., Vignes-Tourneret, F. On a Conjecture of Gross, Mansour and Tucker. Eur. J. Comb., 2021, 97: 103368
\bibitem{ellis2012twisted} Ellis-Monaghan, J.A., Moffatt, I. Twisted Duality for Embedded Graphs. Trans. Am. Math. Soc., 2012, 364: 1529--1569
\bibitem{ellis2013graphs} Ellis-Monaghan, J.A., Moffatt, I. Graphs on Surfaces. Springer, New York, 2013
\bibitem{ellis2011recipe} Ellis-Monaghan, J.A., Sarmiento, I. A Recipe Theorem for the Topological Tutte Polynomial of Bollobás and Riordan. Eur. J. Comb., 2011, 32: 782--794
\bibitem{gross2020partial} Gross, J.L., Mansour, T., Tucker, T.W. Partial Duality for Ribbon Graphs, I: Distributions. Eur. J. Comb., 2020, 86: 103084
\bibitem{gross2021partial} Gross, J.L., Mansour, T., Tucker, T.W. Partial Duality for Ribbon Graphs, II: Partial-Twuality Polynomials and Monodromy Computations. Eur. J. Comb., 2021, 95: 103329
\bibitem{krajewski2011topological} Krajewski, T., Rivasseau, V., Tanasa, A., Wang, Z. Topological Graph Polynomials and Quantum Field Theory, Part II: Mehler Kernel Theories. Ann. Henri Poincaré, 2011, 12: 1--63
\bibitem{wilson1979operators} Wilson, S. Operators over Regular Maps. Pac. J. Math., 1979, 81: 559--568
\bibitem{YDM2026Vertex} Yan, Q., Deng, Q., Metsidik, M. Introducing a Vertex Polynomial Invariant for Embedded Graphs. Discrete Math., 2026, 349: 115111
\bibitem{yan2021counterexamples} Yan, Q., Jin, X. Counterexamples to a Conjecture by Gross, Mansour and Tucker on Partial-Dual Genus Polynomials of Ribbon Graphs. Eur. J. Comb., 2021, 93: 103285
\bibitem{yan2022partial} Yan, Q., Jin, X. Partial-Dual Genus Polynomials and Signed Intersection Graphs. Forum Math. Sigma, 2022, 10: e69
\bibitem{yan2022twist} Yan, Q., Jin, X. Twist Polynomials of Delta-Matroids. Adv. Appl. Math., 2022, 139: 102363
\bibitem{yan2024partial} Yan, Q., Jin, X. Partial-Twuality Polynomials of Delta-Matroids. Adv. Appl. Math., 2024, 153: 102623
\bibitem{yuschak2024delta} Yuschak, D. Delta-Matroids Whose Twist Polynomials Are Monomials. Eur. J. Comb., 2024, 118: 103925

\end{thebibliography}
\end{document}